\documentclass[runningheads,11pt]{llncs}

\usepackage[margin=1in]{geometry}
\usepackage[T1]{fontenc}
\usepackage{graphicx}
\newenvironment{proofsketch}{\paragraph{Proof sketch}}{}

\usepackage{mathtools}
\usepackage{amsfonts}
\usepackage{booktabs}
\usepackage{placeins}
\usepackage[hypertexnames=false,colorlinks=true,breaklinks=true,bookmarks=true,urlcolor=blue,citecolor=blue,linkcolor=blue,bookmarksopen=false,draft=false]{hyperref}

\usepackage{subcaption}
\DeclareMathOperator{\Diag}{Diag}
\DeclareMathOperator{\diag}{diag}
\DeclareMathOperator{\ldet}{ldet}

\DeclareMathOperator{\rank}{rank}

\newtheorem{thm}{Theorem}{\bf}{\rm}
{\bf}{\rm}
\newtheorem{lem}[thm]{Lemma}{\bf}{\rm}
{\bf}{\rm}
\newtheorem{defn}[thm]{Definition}{\bf}{\rm}
{\bf}{\rm}
{\bf}{\rm}

\begin{document}
\title{Generalized Scaling for the Constrained\\ Maximum-Entropy Sampling Problem}
\titlerunning{Generalized Scaling for MESP}
 %
 \author{Zhongzhu Chen\inst{1}
 \and
 Marcia Fampa\inst{2}
 \and
 Jon Lee\inst{1}
 }
 \authorrunning{Chen, Fampa \& Lee}
 %
 \institute{University of Michigan, Ann Arbor, Michigan, USA
 \email{\{zhongzhc,jonxlee\}@umich.edu}\\
  \and
 Universidade Federal do Rio de Janeiro, Brasil\\
 \email{fampa@cos.ufrj.br}}

\maketitle              
\begin{abstract}
The best techniques for the constrained maximum-entropy sampling problem, a discrete-optimization problem arising in the design of experiments, are via a variety of concave continuous relaxations of the objective function. A standard bound-enhancement technique in this context is \emph{scaling}. We extend this technique to \emph{generalized scaling},
we give mathematical results aimed at supporting algorithmic methods for computing optimal generalized scalings, and we give computational results demonstrating the usefulness of generalized scaling on benchmark problem instances.

\keywords{nonlinear $0/1$-optimization  \and convex relaxation \and maximum-entropy sampling}
\end{abstract}
\section{Introduction}\label{sec:intro}

Let $C$ be a symmetric positive semidefinite matrix with rows/columns
indexed from $N:=\{1,2,\ldots,n\}$, with $n >1$.
Let $A\in\mathbb{R}^{m\times n}$ and $b\in\mathbb{R}^m$.
For $0< s < n$,
we define the \emph{constrained maximum-entropy sampling problem}
\begin{equation}\tag{CMESP}\label{CMESP}
\begin{array}{ll}
z:=&\max \left\{\strut\textstyle  \ldet C[S(x),S(x)] ~:~
\mathbf{e}^\top x =s,~ x\in\{0,1\}^n,~Ax\leq b\strut \right\},\\
\end{array}
\end{equation}
where $S(x)$ is the support of $x\in\{0,1\}^n$,
 $C[S,S]$ is the principal submatrix of $C$ indexed by $S$, and
$\ldet$ is the natural logarithm of the determinant.

We refer to \hypertarget{MESP}{MESP} when there are no
constraints   $Ax\leq b$,
which was introduced in the ``design of experiments'' literature by \cite{SW}.
\hyperlink{MESP}{MESP} corresponds to the fundamental problem of choosing an $s$-subvector of
a Gaussian random $n$-vector, so as to maximize the ``differential entropy'' (see \cite{Shannon}).
\hyperlink{MESP}{MESP} has been applied extensively in the field of
environmental monitoring; see \cite[Chapter 4]{FLbook}, and the many references therein.
Important for applications, the constraints $Ax\leq b$ of \ref{CMESP}
can model budget limitations, geographical considerations, and logical dependencies, for example.
We assume $r:=\rank(C)\geq s$, so that
\hyperlink{MESP}{MESP}  always has a feasible solution with finite objective value.

\ref{CMESP} serves as a nice example of a ``non-factorable''
mixed-integer nonlinear program.
When $C$ is a diagonal matrix,
\ref{CMESP} reduces to a general cardinality-constrained
binary linear program.  \cite{AlThaniLeeTridiag,AlThaniLeeTridiag_journal} established that when $C$ is tridiagonal (or even when the support graph of $C$ is a spider with a bounded number of legs), \hyperlink{MESP}{MESP} is then polynomially solvable
by dynamic programming.

\cite{KLQ} established that \hyperlink{MESP}{MESP} is NP-hard
and introduced a novel B\&B (branch-and-bound) approach
based on a spectral bound. \cite{LeeConstrained}
extended the spectral approach to \ref{CMESP}.
\cite{AFLW_IPCO}  and \cite{AFLW_Using}
developed a bound employing a novel convex relaxation.
\cite{Anstreicher_BQP_entropy} developed the  ``BQP bound'', using an extended formulation based on the Boolean quadric polytope. \cite{Kurt_linx} introduced the ``linx bound'', based on a clever convex relaxation.
\cite{nikolov2015randomized} gave a novel
 ``factorization bound'' based on a somewhat mysterious convex
 relaxation. This was further developed by \cite{li2020best} and then \cite{Fact}.
 \cite{chen_mixing} gave a methodology for combining multiple
 convex-optimization bounds to give improved bounds.
 All of these convex-optimization based bounds admit variable
 fixing methodology based on convex duality (see \cite{FLbook}, for example).
 Another key idea for deriving bounds is ``complementation''.
If $C$ is invertible,
we have
\[
z=z(C^{-1},n-s,-A,b-A\mathbf{e}) + \ldet C,
\]
\noindent where $z(C^{-1},n-s,-A,b-A\mathbf{e})$ denotes the
optimal value of \ref{CMESP} with $C,s,A,b$ replaced by
$C^{-1},n-s,-A,b-A\mathbf{e}$, respectively.
So we have a
\emph{complementary} \ref{CMESP} problem
 and \emph{complementary} bounds
 (i.e., bounds for the
 complementary problem plus $\ldet C$)
 immediately give us bounds on
 $z$.
 Some upper bounds on $z$ also shift by  $\ldet C$  under complementing,
 in which case there is no additional value in computing the
 complementary bound.
 Details on all of this can be found in
\cite{FLbook}.

\subsubsection{Terminology.}
Throughout, we let  $\Upsilon:=\left(\gamma_1,\gamma_2,\ldots,\gamma_n\right)^\top \in \mathbb{R}_{++}^n$ be a ``scaling vector''.
We refer to our bounds as \emph{g-scaled} (i.e., \emph{generalized scaled}), and when all elements
of $\Upsilon$ are equal, we  say  \emph{o-scaled} (i.e., \emph{ordinary scaled}).
If all elements of $\Upsilon$ are equal to 1,
we say \emph{un-scaled}.

\subsubsection{Organization and contributions.}
In \S\ref{sec:bqp}, we introduce the g-scaled
BQP bound and establish its convexity in the log of the scaling vector,
generalizing an important and practically-useful result (see \cite[Thm. 11]{chen_mixing}).
In \S\ref{sec:linx}, we introduce the g-scaled
linx bound and establish its convexity in the log of the scaling vector,
generalizing another very important and practically-useful
result for o-scaling (see \cite[Thm. 18]{chen_mixing}).
These convexity results are key
for the tractability of globally optimizing the scaling, something that we
do not have for more general bound ``masking'' (see \cite{AnstreicherLee_Masked,BurerLee}).
In \S\ref{sec:fact}, we introduce the g-scaled
factorization bound, and we establish that g-scaling can significantly improve the factorization bound for \ref{CMESP}, while the o-scaling cannot help it (see \cite[Thm. 2.1]{Fact}). We are also able to prove that
for \hyperlink{MESP}{MESP}, the all-ones vector is a stationary
point for the bound as a function of the scaling vector. Therefore,
g-scaling is unlikely to be helpful for \hyperlink{MESP}{MESP}, similar to o-scaling.
In \S\ref{sec:num}, we present results of computational experiments, demonstrating the improvements on upper bounds
and on the number of variables that can be fixed (using convex duality)
due to g-scaling.
In \S\ref{sec:conc}, we make some brief concluding remarks.
In \S\ref{sec:proofs}, we provide some proof sketches.

\subsubsection{Notation.}
$\Diag(x)\in\mathbb{R}^{n\times n}$ makes a diagonal matrix from
$x\in\mathbb{R}^n$. $\diag(X)\in\mathbb{R}^n$ extracts the diagonal of $X\in\mathbb{R}^{n\times n}$.
We let $\mathbb{S}^n_+$ (resp., $\mathbb{S}^n_{++}$) be the set of
 positive semidefinite (resp., definite) symmetric matrices of order $n$.
 We let $\lambda_\ell(M)$ be the $\ell$-th greatest eigenvalue of
$M\in \mathbb{S}^n_+$\thinspace. We denote by $\mathbf{e}$ an all-ones vector.
For matrices $A$ and $B$ with the same shape,
$A\circ B $ is the Hadamard (i.e., element-wise) product. We denote natural logarithm by $\log$,
and apply it component-wise to vectors.


\section{BQP bound}\label{sec:bqp}

We define  the convex set
\begin{align*}
&P(n,s):=\left\{
(x,X)\in\mathbb{R}^n\times \mathbb{S}^n ~:~
 X -xx^\top\succeq 0,~ \diag(X)=x,~
 \mathbf{e}^\top x=s,~ X\mathbf{e}=sx
\right\}.
\end{align*}
For $\Upsilon \in \mathbb{R}_{++}^n$,  $x\in[0,1]^n$ and $X\in \mathbb{S}_+^n$, we define
   \begin{align*}
&	f_{{\tiny\mbox{BQP}}}(x,X;\Upsilon):=\textstyle \ldet \left( \left(\Diag(\Upsilon)C\Diag(\Upsilon)\right)\circ X + \Diag(\mathbf{e}-x)\right) - 2\sum_{i=1}^n  x_i\log \gamma_i 
	\end{align*}
	
\noindent and the \emph{g-scaled BQP bound}
	\[
	z_{{\tiny\mbox{BQP}}}(\Upsilon):=\max \left\{f_{{\tiny\mbox{BQP}}}(x,X;\Upsilon) ~:~(x,X)\in P(n,S),~  Ax\leq b \right\}.\tag{BQP}\label{BQP}
	\]
Note that we can interpret this bound as applying the
un-scaled BQP bound to the symmetrically-scaled
matrix $\Diag(\Upsilon)C\Diag(\Upsilon)$, and then
correcting by $- 2\sum_{i=1}^n  x_i\log \gamma_i$~.

\begin{thm}\label{thm:bqp}
\phantom{.}
\vspace{-10pt}

\begin{itemize}
\item[\ref{thm:bqp}.i.] \label{bqp.i} $z\leq z_{{\tiny\mbox{BQP}}}$~;
\item[\ref{thm:bqp}.ii.] For all $\Upsilon \in \mathbb{R}_{++}^n$, $f_{{\tiny\mbox{BQP}}}(x,X;\Upsilon)$  is concave on the feasible region of \ref{BQP};
\item[\ref{thm:bqp}.iii.] $z_{{\tiny\mbox{BQP}}}(\Upsilon)$ is convex in $\log \Upsilon$.
\end{itemize}
\end{thm}

\noindent The \ref{BQP} bound was first analyzed and developed in \cite{Anstreicher_BQP_entropy}, establishing Thm. \ref{thm:bqp}.\emph{i} for $\Upsilon=\mathbf{e}$. Thm. \ref{thm:bqp}.\emph{ii}
is a result of \cite{Anstreicher_BQP_entropy}, with details
filled in by \cite{FLbook}. Thm. \ref{thm:bqp}.\emph{iii}
significantly generalizes a result of \cite{chen_mixing}, where it is
established only for o-scaling: i.e., on $\{ \Upsilon= \gamma \mathbf{e} ~:~ \gamma\in \mathbb{R}_{++}\}$. The proof of Thm. \ref{thm:bqp}.\emph{iii} requires new ideas (see the proof sketch in the Appendix). Additionally, the result is quite important as it enables the use of readily available quasi-newton methods (like BFGS) for finding the globally optimal g-scaling for the \ref{BQP} bound.

\section{linx bound}\label{sec:linx}

For  $\Upsilon \in \mathbb{R}_{++}^n$ and  $x\in[0,1]^n$,
 we define
\[
f_{{\tiny\mbox{linx}}}(x;\Upsilon) :=
\textstyle{\frac{1}{2}}\left(\ldet \left( \Diag(\Upsilon) C\Diag(x)C \Diag(\Upsilon)+\Diag(\mathbf{e}-x) \right)\right) -\sum_{i=1}^n x_i\log \gamma_i
\]
and the \emph{g-scaled linx  bound}  
\begin{equation}\tag{linx}\label{linx}
\begin{array}{ll}
	z_{{\tiny\mbox{linx}}}(\Upsilon):=&\max \left\{
	\strut
	f_{{\tiny\mbox{linx}}}(x;\Upsilon)
	~:~
 \mathbf{e}^{\top}x=s,~0\leq x\leq \mathbf{e},~ Ax\leq b\strut \right\}.
\end{array}
\end{equation}
Note that we \emph{cannot} interpret this bound as applying the
un-scaled linx bound to the row-scaled
matrix $\Diag(\Upsilon)C$, because we would lose symmetry.

\vbox{
\begin{thm}\label{thm:linx}
\phantom{.}
\vspace{-10pt}

\begin{itemize}
\item[\ref{thm:linx}.i.] $z\leq z_{{\tiny\mbox{linx}}}$~;
\item[\ref{thm:linx}.ii.] For all $\Upsilon \in \mathbb{R}_{++}^n$, $f_{{\tiny\mbox{linx}}}(x;\Upsilon)$  is concave on the feasible region of \ref{linx};
\item[\ref{thm:linx}.iii.] $z_{{\tiny\mbox{linx}}}(\Upsilon)$ is convex in $\log \Upsilon$.
\end{itemize}
\end{thm}
}

\noindent The \ref{linx} bound was first analyzed and developed in \cite{Kurt_linx}, establishing Thm. \ref{thm:linx}.\emph{i} for $\Upsilon=\mathbf{e}$. Thm. \ref{thm:linx}.\emph{ii}
is a result of \cite{Kurt_linx}, with details
filled in by \cite{FLbook}. Thm. \ref{thm:linx}.\emph{iii}
generalizes a result of \cite{chen_mixing}, where it is
established only for o-scaling: i.e., on $\{ \Upsilon= \gamma \mathbf{e} ~:~ \gamma\in \mathbb{R}_{++}\}$. The proof of Thm. \ref{thm:linx}.\emph{iii} requires new ideas (see the proof sketch in the Appendix). Additionally, the result is quite important as it enables the use of readily available quasi-newton methods (like BFGS) for finding the globally optimal g-scaling for the \ref{linx} bound.

\section{Factorization bound}\label{sec:fact}


\begin{lem}\label{Nik}(see \cite[Lem. 14]{nikolov2015randomized})
 Let $\lambda\in\mathbb{R}_+^k$ with $\lambda_1\geq \lambda_2\geq \cdots\geq \lambda_k$~, and let $0<s\leq k$. There exists a unique integer $\iota$, with $0\leq \iota< s$, such that
 $
 \lambda_{\iota}>\frac{1}{s-\iota}\sum_{\ell=\iota+1}^k \lambda_{\ell}\geq \lambda_{\iota +1}
 $~,
 with the convention $\lambda_0=+\infty$.
\end{lem}
Now, suppose that  $\lambda\in\mathbb{R}^k_+$ with
$\lambda_1\geq\lambda_2\geq\cdots\geq\lambda_k$. Given an integer $s$ with $0<s\leq k$,
let $\iota$ be the unique integer defined by Lem. \ref{Nik}. We define
$
\phi_s(\lambda):=\sum_{\ell=1}^{\iota} \log \lambda_\ell + (s - \iota)\log\left(\frac{1}{s-{\iota}} \sum_{\ell=\iota+1}^{k}
\lambda_\ell\right)
$.
Next, for $X\in\mathbb{S}_{+}^k$~, we define
$\Gamma_s(X):=  \phi_s(\lambda_1(X),\ldots,\lambda_k(X))$.

Suppose that the rank of $C$ is $r\geq s$.
Then we factorize $C=FF^\top$,
with $F\in \mathbb{R}^{n\times k}$, for some $k$ satisfying $r\le k \le n$.
Now, for
$\Upsilon \in \mathbb{R}_{++}^n$ and
$x\in [0,1]^n$, we define
$F_{{\tiny\mbox{DDFact}}}(x;\Upsilon):
= \sum_{i=1}^n \gamma_i x_i F_{i\cdot}^\top F_{i\cdot}$~.

\noindent Finally, we define
$
f_{{\tiny\mbox{DDFact}}}(x;\Upsilon):= \strut \Gamma_s(F_{{\tiny\mbox{DDFact}}}(x;\Upsilon)) -\sum_{i=1}^n x_i \log \gamma_i
$
and the \emph{g-scaled factorization bound}
\[
\begin{array}{ll}
	z_{{\tiny\mbox{DDFact}}}(\Upsilon):= 	&\max \left\{  f_{{\tiny\mbox{DDFact}}}(x;\Upsilon)~:~  \mathbf{e}^\top x=s,~
 0\leq x\leq \mathbf{e}, ~ Ax\leq b\strut \right\}.
\end{array}
\tag{DDFact}\label{DDFact}
\]
Noticing that $F_{{\tiny\mbox{DDFact}}}(x;\Upsilon):
=F^\top \Diag\left(\sqrt{\Upsilon}\right) \Diag(x) \Diag\left(\sqrt{\Upsilon}\right)F$,
we can interpret this bound as applying the
un-scaled \ref{DDFact} bound to the symmetrically-scaled
matrix $\Diag\left(\sqrt{\Upsilon}\right)C\Diag\left(\sqrt{\Upsilon}\right)$, and then
correcting by $- \sum_{i=1}^n  x_i\log \gamma_i$~.


\begin{defn}\label{def:fact}
For any $x$ feasible to \ref{DDFact}, suppose the eigenvalues of $F_{{\tiny\mbox{DDFact}}}(x;\Upsilon)$ are $\lambda_1\ge \cdots\ge \lambda_r > \lambda_{r+1}=\cdots=\lambda_{k}
=0$~, where $r\in [s,k]$ and $F_{{\tiny\mbox{DDFact}}}(x;\Upsilon) =  Q \Diag( \lambda ) Q $ with an orthonormal matrix $ Q $. Define $ \beta ( \lambda ):=(\beta_1,\beta_2,\ldots,\beta_{k})^\top$ such that
\begin{align*}
   \beta_{i}:=\frac{1}{\lambda_{i}},~ \forall~ i \in[1, \iota],~ \beta_{i}:=\frac{s-\iota}{\sum_{i \in[\iota+1, k]} \lambda_{i}},~ \forall~ i \in[\iota+1, k],
\end{align*}
where $\iota$ is the unique integer in Lemma \ref{Nik}.
\end{defn}

\begin{thm}\label{thm:fact}
\phantom{.}
\vspace{-10pt}

\begin{itemize}
\item[\ref{thm:fact}.i.] $z\leq z_{{\tiny\mbox{DDFact}}}$~;
\item[\ref{thm:fact}.ii.] For all $\Upsilon  \in \mathbb{R}_{++}^n$, $f_{{\tiny\mbox{DDFact}}}(x;\Upsilon)$ is concave on the feasible region of \ref{DDFact};
\item[\ref{thm:fact}.iii.] For all $\Upsilon  \in \mathbb{R}_{++}^n$ and $x\geq 0$ in the domain of $f_{{\tiny\mbox{DDFact}}} (x;\Upsilon)$, let $T(x;\Upsilon):=$\break $ \diag\left(F_{{\tiny\mbox{DDFact}}}(x;\Upsilon) Q \Diag\left( \beta ( \lambda )\right) Q ^\top F_{{\tiny\mbox{DDFact}}}(x;\Upsilon)^\top\right) - \log\Upsilon $ where $Q, \beta(\lambda)$ are defined in \ref{def:fact}, then
\begin{align*}
    \lim_{
    \|\hat x-x\|\rightarrow 0 ~:~\atop
    {\hat x \geq 0 \text{ is in the domain}
    \atop
    {\text{ of } f_{{\tiny\mbox{DDFact}}} (x;\Upsilon)}}
    } \frac{\left|f_{{\tiny\mbox{DDFact}}}(\hat x;\Upsilon)-f_{{\tiny\mbox{DDFact}}}(x;\Upsilon)-T(x;\Upsilon)^\top (\hat x-x)\right|}{\|\hat x-x\|} = 0.
\end{align*}

\item[\ref{thm:fact}.iv.] For all $x$ feasible in the domain of $f_{{\tiny\mbox{DDFact}}} (x;\Upsilon)$, $f_{{\tiny\mbox{DDFact}}}(x;\Upsilon)$ is differentiable in $\Upsilon$ at all $\Upsilon  \in \mathbb{R}_{++}^n$. In particular, for \hyperlink{MESP}{MESP}, let $x^*$ to be one optimal solution to \ref{DDFact}, then we have
\begin{align*}
    \left.\frac{\partial f_{{\tiny\mbox{DDFact}}} (x^*;\Upsilon)}{ \partial \Upsilon}\right|_{\Upsilon=\mathbf{e}} = 0.
\end{align*}
\end{itemize}
\end{thm}

\noindent The \ref{DDFact} bound was first analyzed and developed in \cite{nikolov2015randomized},
establishing Thm. \ref{thm:fact}.\emph{i} for $\Upsilon=\mathbf{e}$,
and developed further in \cite{li2020best}.
We note that the o-scaled factorization bound for \ref{CMESP} is
invariant under the scale factor (see \cite{Fact}), so the
use of any type of scaling in the context of the \ref{DDFact} bound is completely new.
Thms. \ref{thm:fact}.\emph{iii}-\emph{iv} are the first differentiablity results of any type for the
\ref{DDFact} bound. The proof methods (sketched in the Appendix) are quite technical and novel. Furthermore, they explain the success of our quasi-newton based methods for calculating
optimal g-scalings for the \ref{DDFact} bound, not anticipated by previous works which exposed only subgradients connected to   \ref{DDFact}.
As we will see in \S\ref{sec:num}, g-scaling can improve
the \ref{DDFact} bound for  \ref{CMESP}.
These observations and Thm. \ref{thm:fact}.\emph{iv} leave open
the interesting question of whether g-scaling can help
the \ref{DDFact} bound for \hyperlink{MESP}{MESP};
we can interpret Thm. \ref{thm:fact}.\emph{iv} as a partial result toward a negative answer.

\section{Numerical results}\label{sec:num}

We experimented on benchmark instances of \hyperlink{MESP}{MESP}, using three covariance matrices
 that have been extensively used in the literature, with $n=63,90,124$ (see, e.g., \cite{KLQ,LeeConstrained,AFLW_Using,Anstreicher_BQP_entropy,Kurt_linx}).
 For testing \ref{CMESP}, we included five side constraints $a_i^\top x\leq b_i$, for $i=1,\ldots,5$, in \hyperlink{MESP}{MESP}. As there is no benchmark data for the side constraints, we have generated them randomly.  For each $n$, the  left-hand side of constraint $i$ is given by a uniformly-distributed random vector $a_i$ with integer components between $-2$ and $2$. The right-hand side of the constraints was selected so that, for
every $s$ considered in the experiment, the best known solution
of the instance of \hyperlink{MESP}{MESP} is violated by
at least one constraint.



For each $n$,
we considered  instances of \hyperlink{MESP}{MESP} and \ref{CMESP} with
a wide range of $s$.
We ran our experiments under Windows, on an Intel Xeon E5-2667 v4 @ 3.20 GHz processor equipped with 8 physical cores (16 virtual cores) and 128 GB of RAM.
We implemented our code in \texttt{Matlab} using the solvers
 \texttt{SDPT3} v. 4.0 for \ref{BQP}, and \texttt{Knitro} v. 12.4
for \ref{linx} and  \ref{DDFact}, and optimizing scaling vectors $\Upsilon$ using a BFGS algorithm, and the o-scaling parameters $\gamma$ using the Newton's method. Besides solving the relaxations to get upper bounds for our test instances of
\hyperlink{MESP}{MESP} and \ref{CMESP}, we  compute lower bounds with a heuristic of \cite[Sec. 4]{LeeConstrained} and then a local search (see \cite[Sec. 4]{KLQ}).

\begin{figure*}
        \centering
        \begin{subfigure}[b]{0.496\textwidth}
            \centering
            \includegraphics[height=5cm,width=\textwidth]{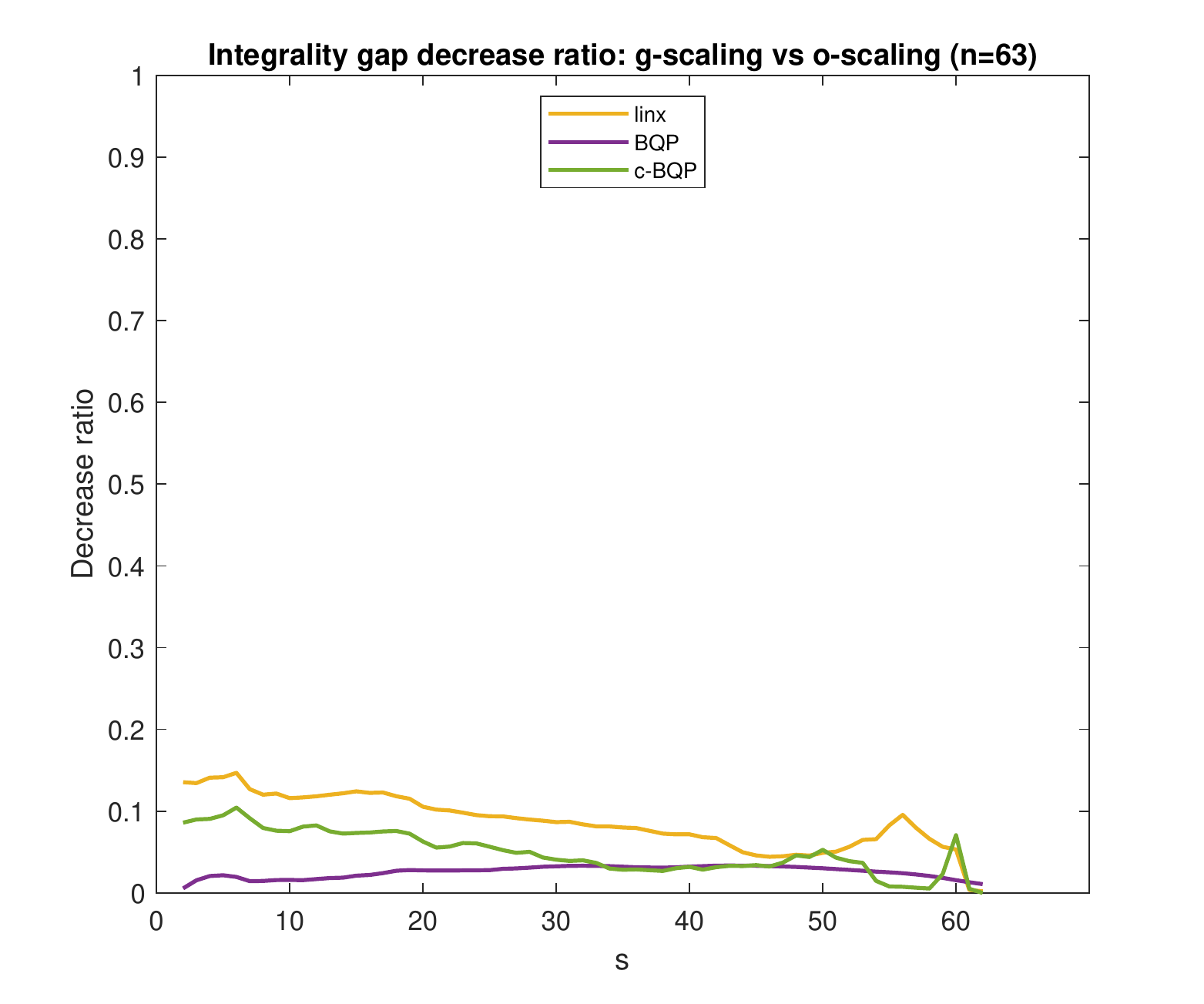}
        \end{subfigure}
        \hfill
        \begin{subfigure}[b]{0.496\textwidth}
            \centering
            \includegraphics[height=5cm,width=\textwidth]{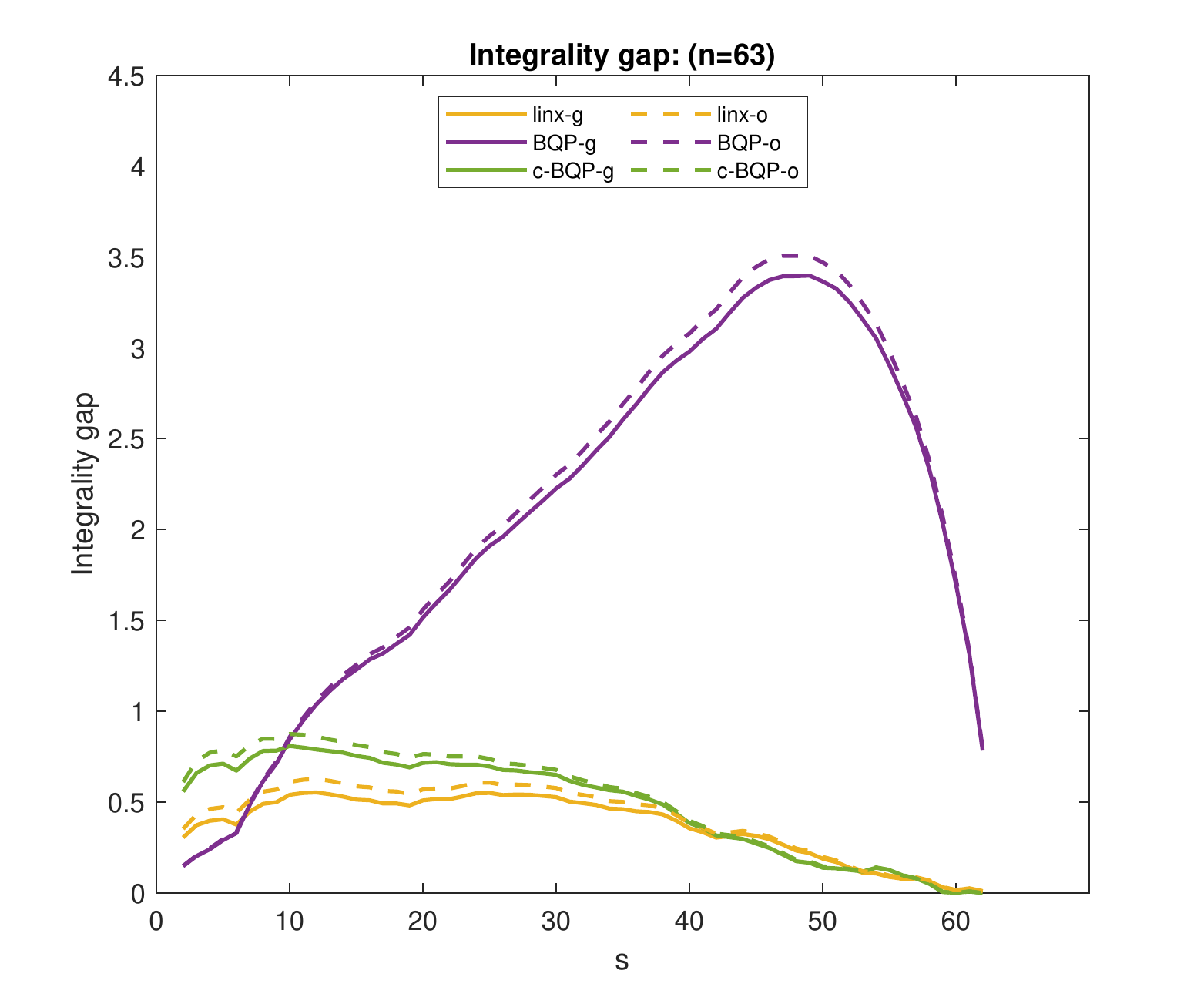}
        \end{subfigure}
        \vskip\baselineskip
        \begin{subfigure}[b]{0.496\textwidth}
            \centering
            \includegraphics[height=5cm,width=\textwidth]{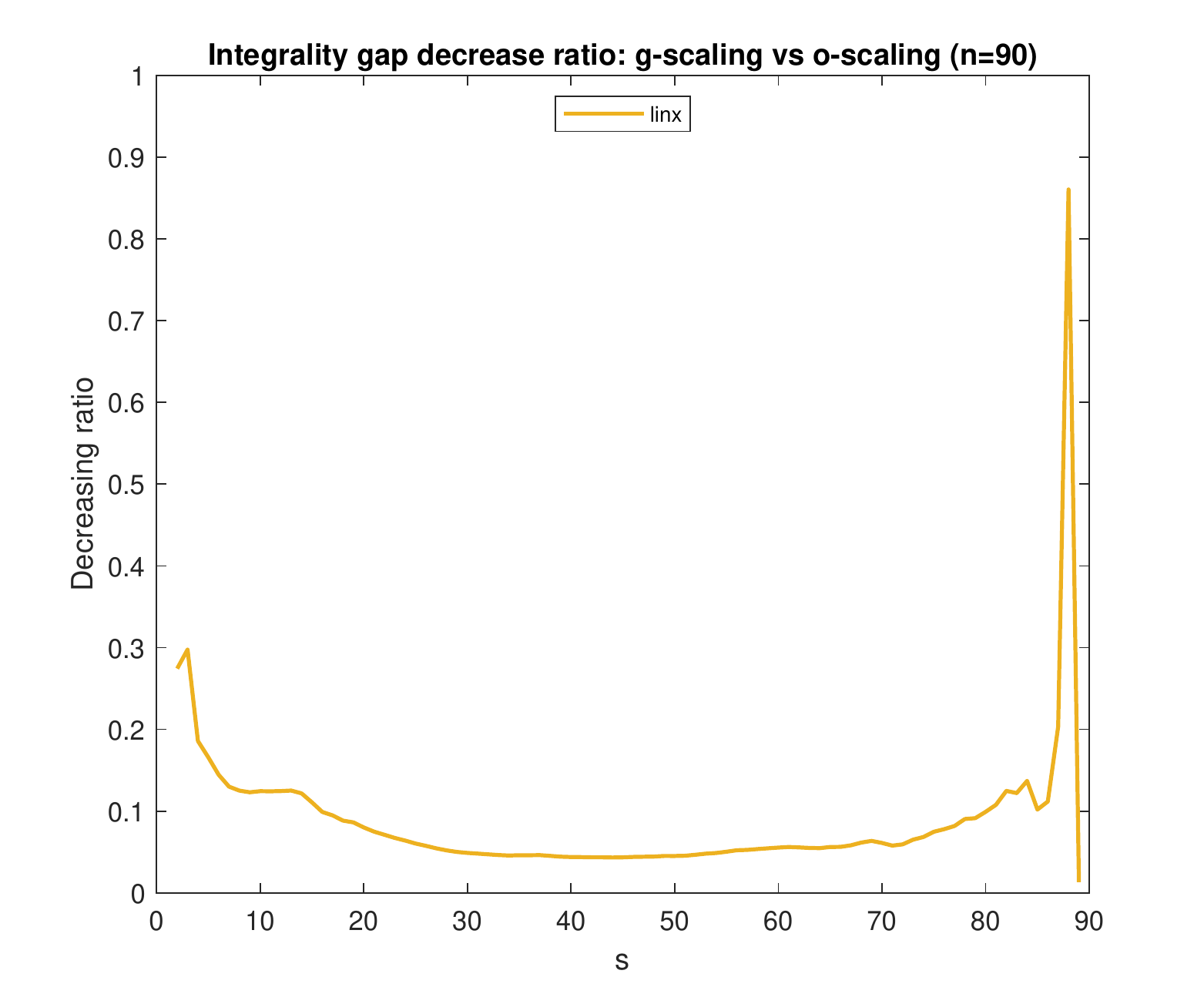}
        \end{subfigure}
        \hfill
        \begin{subfigure}[b]{0.496\textwidth}
            \centering
            \includegraphics[height=5cm,width=\textwidth]{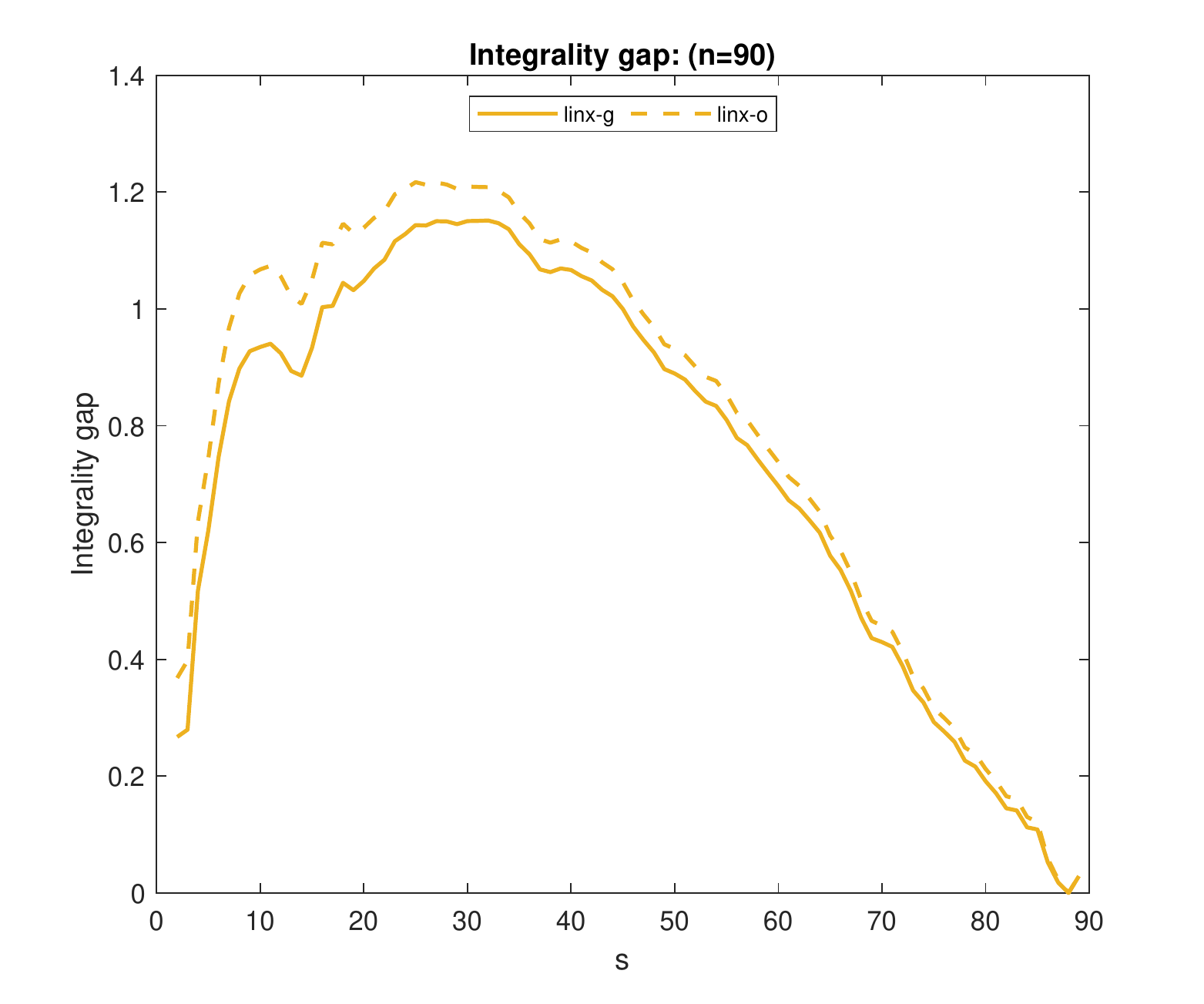}
        \end{subfigure}
        \vskip\baselineskip
        \begin{subfigure}[b]{0.496\textwidth}
            \centering
            \includegraphics[height=5cm,width=\textwidth]{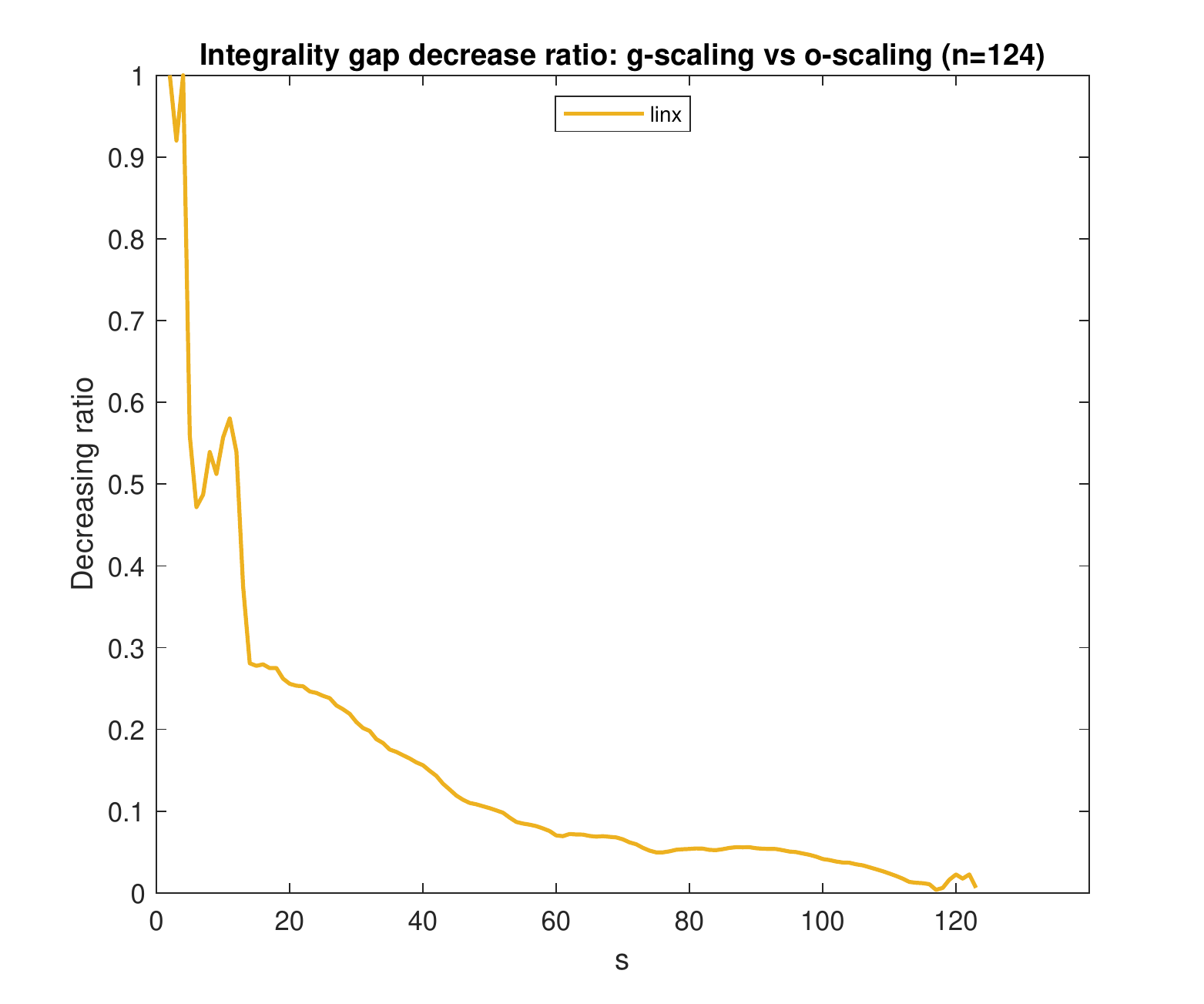}
        \end{subfigure}
        \hfill
        \begin{subfigure}[b]{0.496\textwidth}
            \centering
            \includegraphics[height=5cm,width=\textwidth]{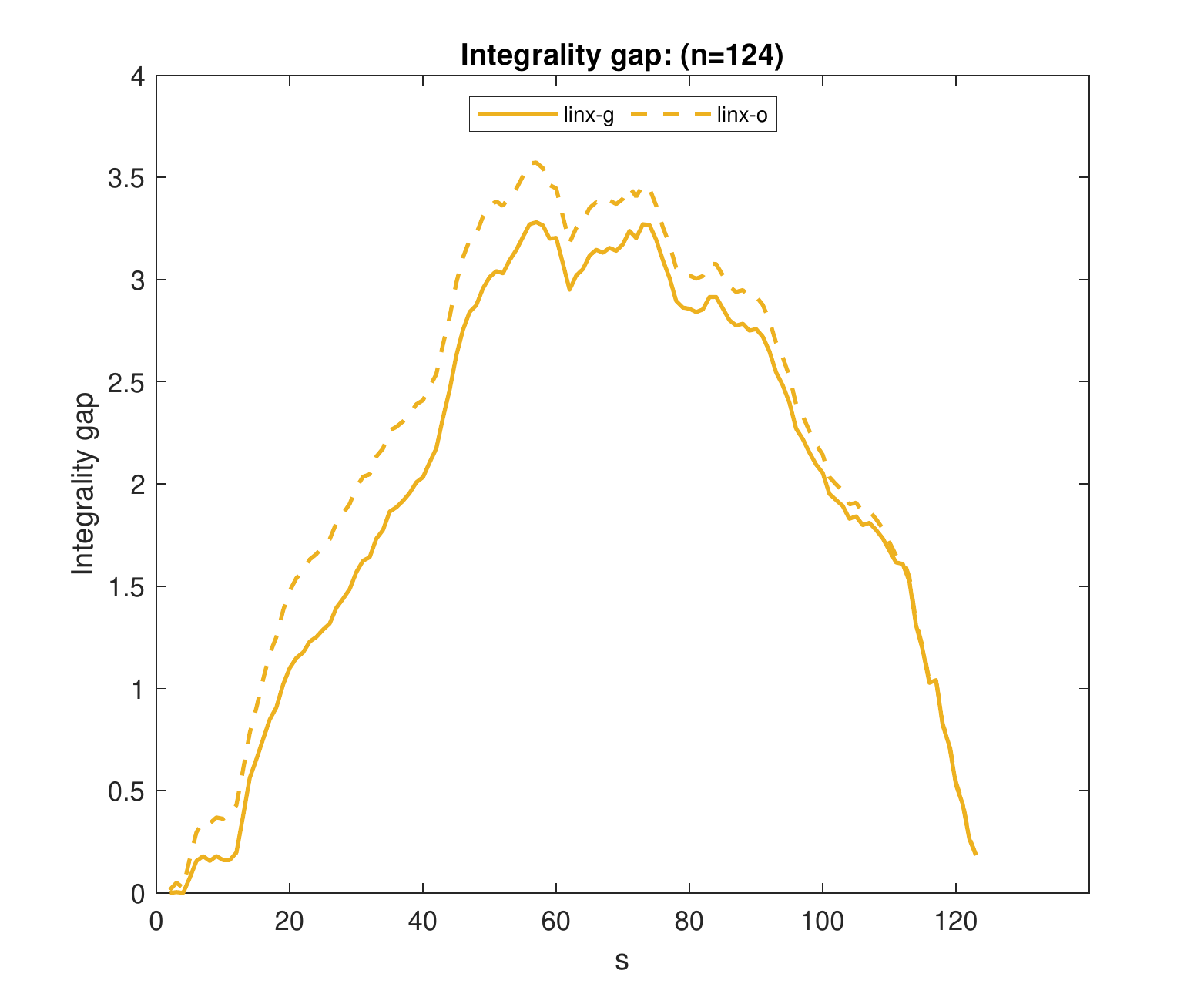}
        \end{subfigure}


        \caption[ ]
        {\small Comparison between g-scaling and o-scaling for \hyperlink{MESP}{MESP}}
        \label{fig:unconst}
    \end{figure*}

In Fig. \ref{fig:unconst}, we show the impact of  g-scaling on the  linx bound for \hyperlink{MESP}{MESP} on the three benchmark covariance matrices. For  the $n=63$ matrix, we also show the impact of  g-scaling on the BQP bound. The DDFact and complementary DDFact bounds are only considered in the experiments for \ref{CMESP}, as the g-scaling methodology was only able to improve these bounds when side constraints were added to \hyperlink{MESP}{MESP}. The plots on the left in Fig. \ref{fig:unconst} present the ``integrality gap decrease ratios'', given by the difference between the integrality gaps using o-scaling and the integrality gaps using g-scaling, divided by the  integrality gaps using o-scaling.
The integrality gaps are given by the difference between the upper bounds computed with the relaxations and
     lower bounds given by heuristic solutions. We see that larger  $n$ leads to larger  maximum ratios. We also see that the g-scaling methodology is effective in reducing all bounds evaluated, especially the linx bound. Even for the most difficult instances, with intermediate values of $s$, we have some improvement on the bounds, which can be  effective in the branch-and-bound context where the bounds would ultimately be applied. The plots on the right in Fig. \ref{fig:unconst} present the integrality gaps, and we see that even when the integrality gaps given by the o-scaling are less than 1,  g-scaling can reduce them.

\begin{figure*}
        \centering
        \begin{subfigure}[b]{0.496\textwidth}
            \centering
            \includegraphics[height=5cm,width=\textwidth]{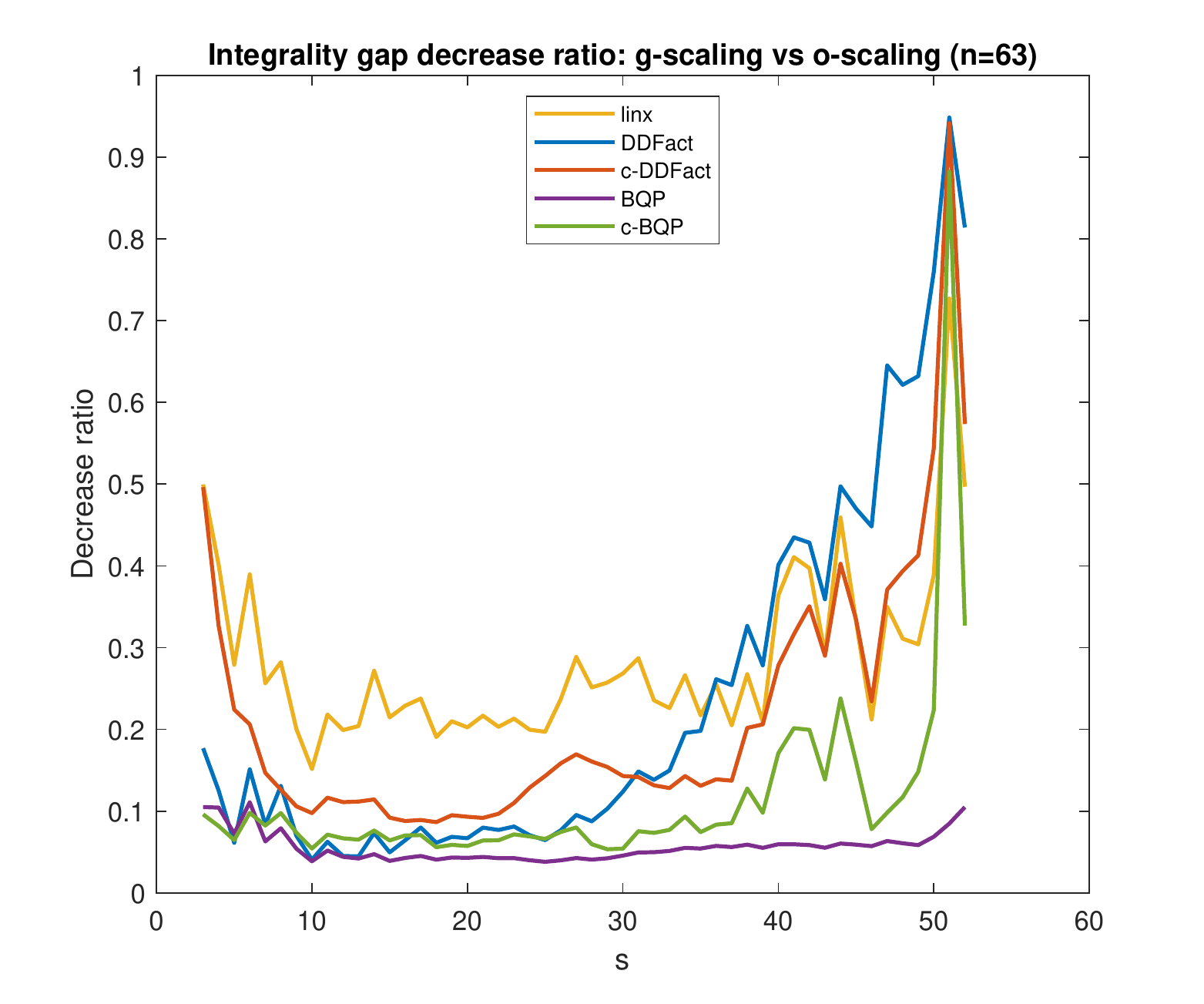}
        \end{subfigure}
        \hfill
        \begin{subfigure}[b]{0.496\textwidth}
            \centering
            \includegraphics[height=5cm,width=\textwidth]{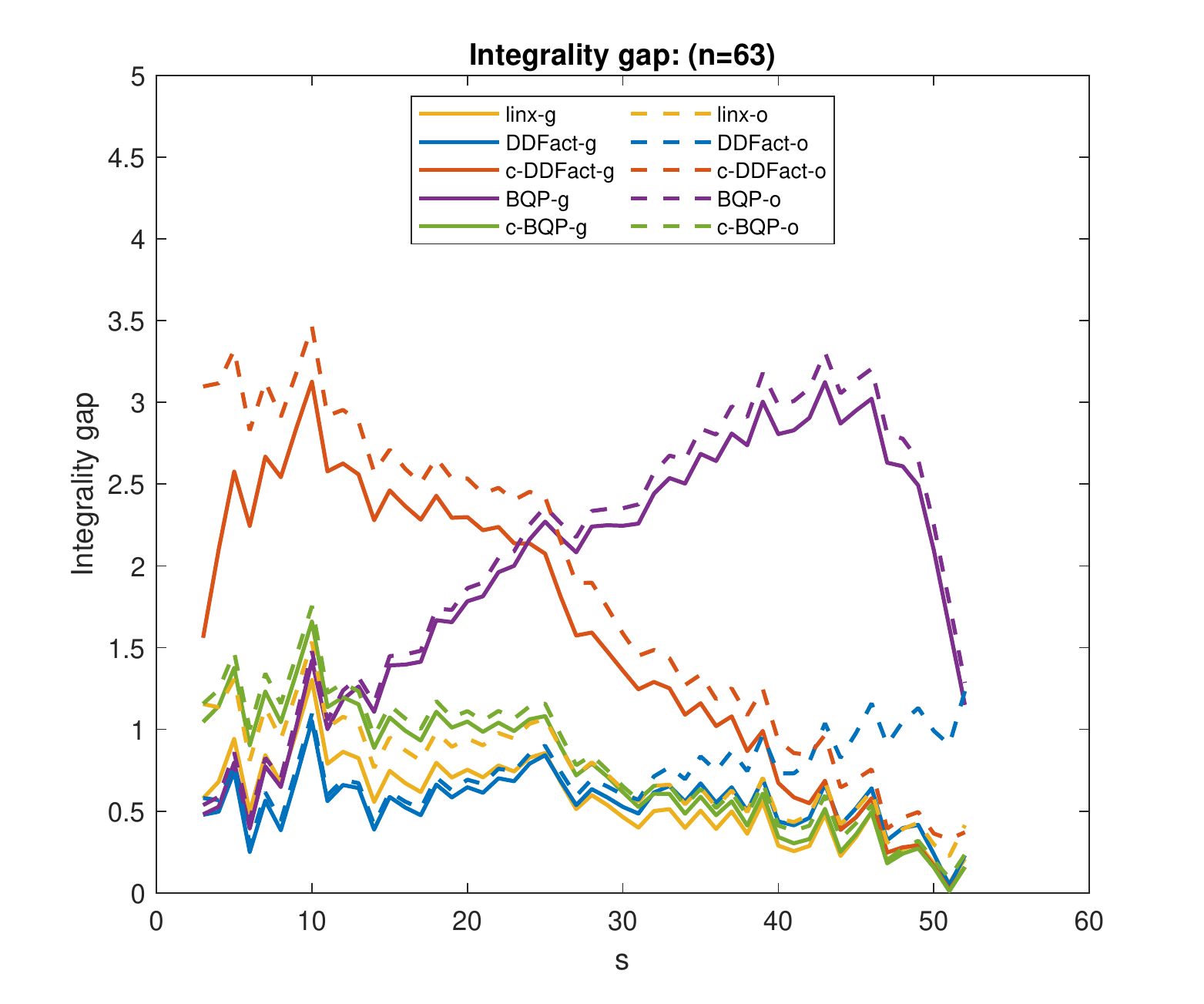}
        \end{subfigure}
        \vskip\baselineskip
        \begin{subfigure}[b]{0.496\textwidth}
            \centering
            \includegraphics[height=5cm,width=\textwidth]{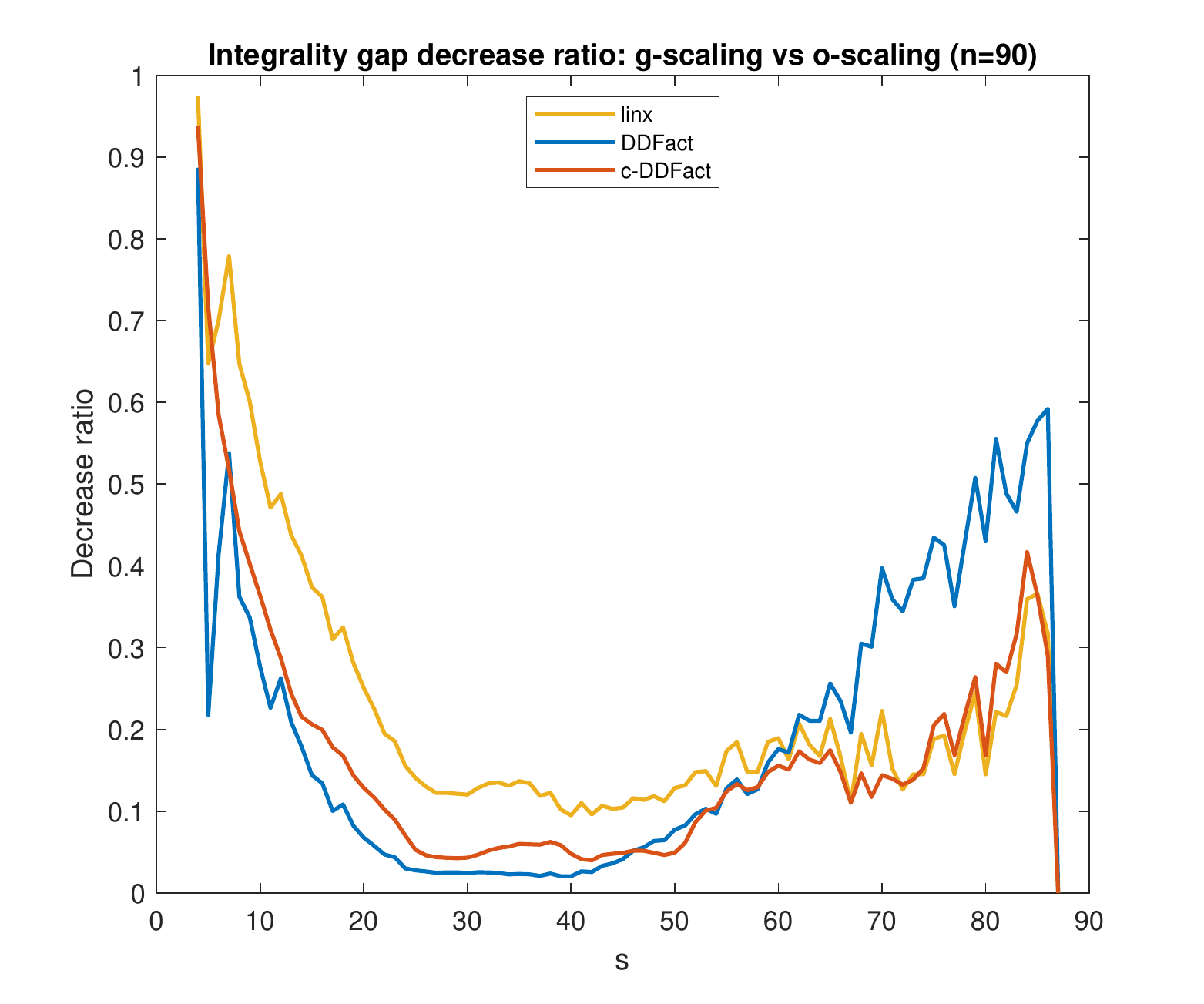}
        \end{subfigure}
        \hfill
        \begin{subfigure}[b]{0.496\textwidth}
            \centering
            \includegraphics[height=5cm,width=\textwidth]{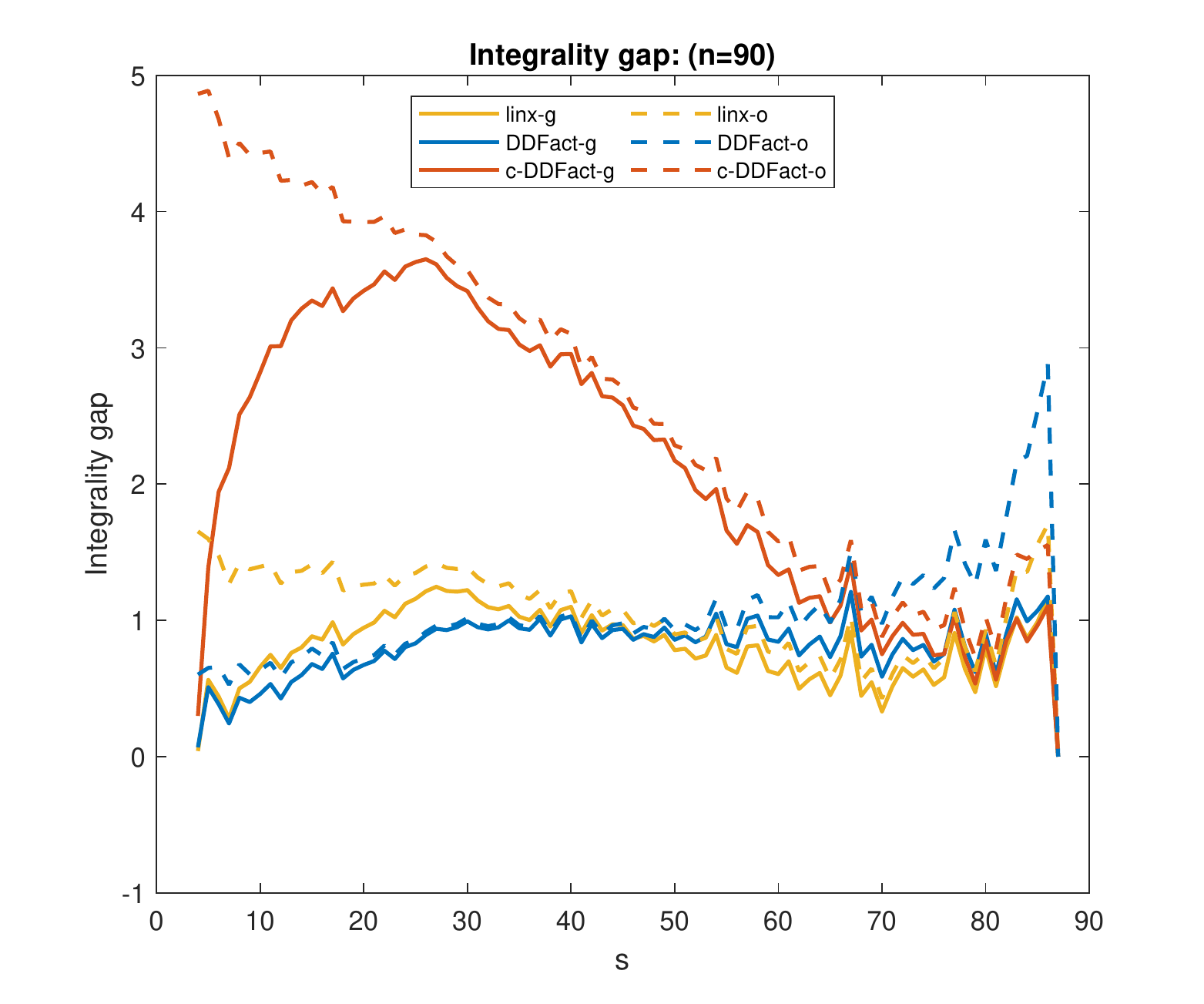}
        \end{subfigure}
        \vskip\baselineskip
        \begin{subfigure}[b]{0.496\textwidth}
            \centering
            \includegraphics[height=5cm,width=\textwidth]{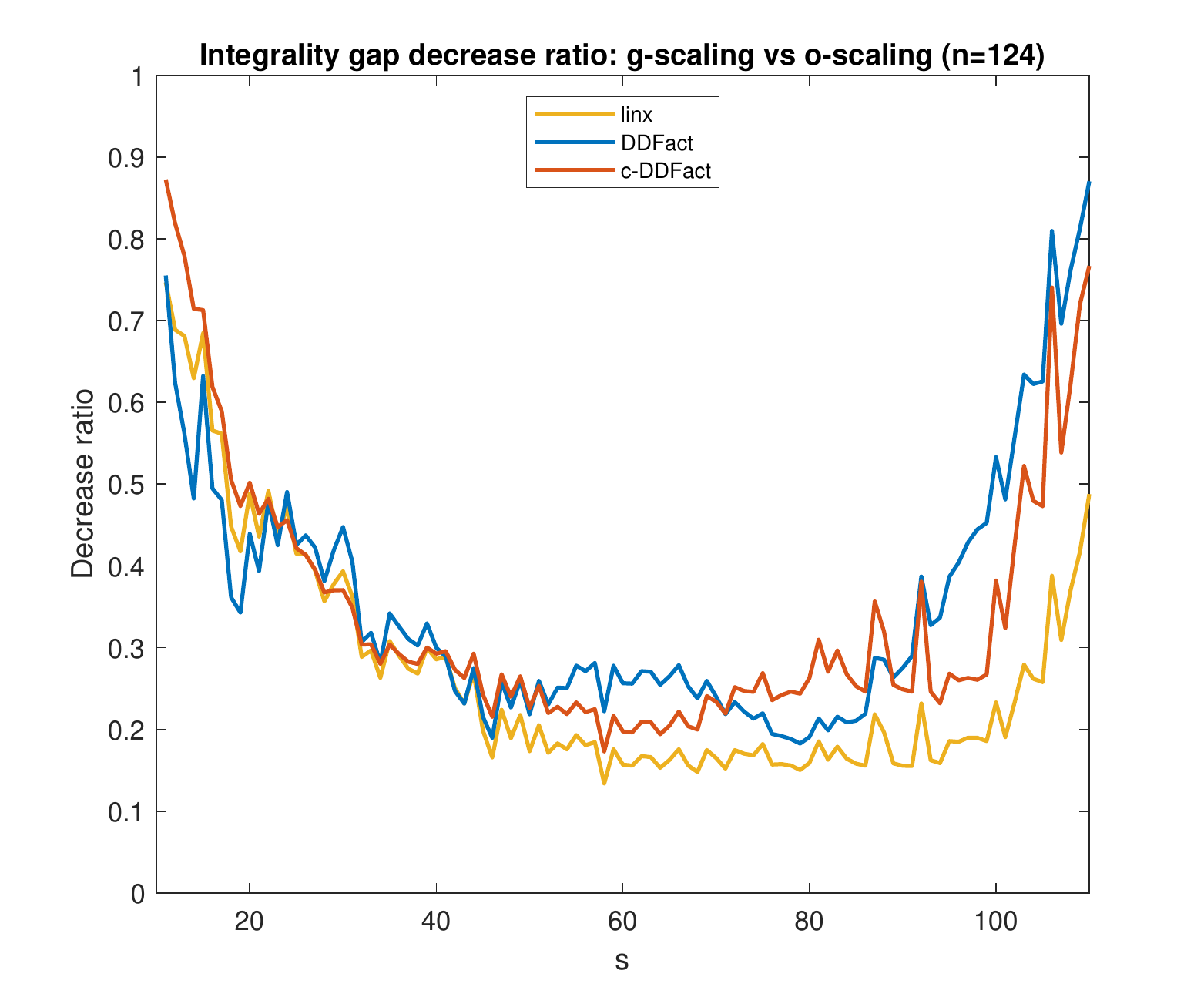}
        \end{subfigure}
        \hfill
        \begin{subfigure}[b]{0.496\textwidth}
            \centering
            \includegraphics[height=5cm,width=\textwidth]{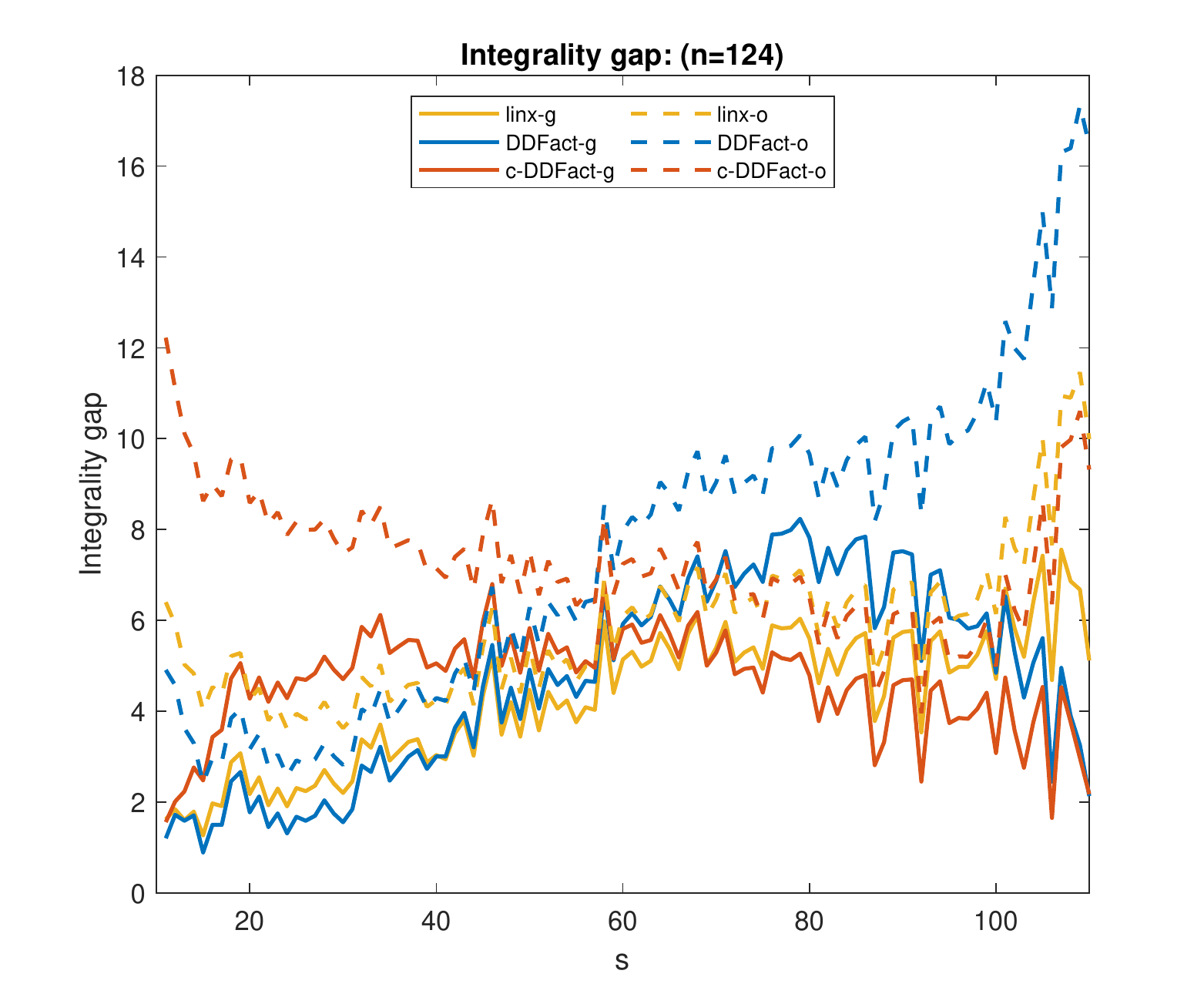}
        \end{subfigure}
        \caption[ ]
        {\small Comparison between g-scaling and o-scaling for \ref{CMESP}}
        \label{fig:const}
    \end{figure*}

In Fig. \ref{fig:const}, we show for \ref{CMESP}, similar results  to the ones shown in Fig. \ref{fig:unconst}, except that now we also present the effect of g-scaling on  the DDFact and the complementary DDFact bounds. We see from the integrality gap decrease ratios that when side constraints are added to  \hyperlink{MESP}{MESP}, the g-scaling is, in general, more effective in reducing the gaps given by o-scaling. We also see that,   it is particularly  effective in reducing the DDFact and complementary DDFact bounds. Especially for the $n=124$ matrix, we see a significant reduction on the gaps given by complementary DDFact and DDFact, for $s$ smaller and greater than  $50$, respectively.

We also investigated how the improvement of g-scaling  over o-scaling for the linx bound can increase the possibility of fixing variables in \hyperlink{MESP}{MESP} and \ref{CMESP}. The methodology for  fixing variables is based on convex duality and  has been applied  since the first convex relaxation  was proposed for these problems in \cite{AFLW_IPCO}. When a lower bound for each problem is available, the dual solution of the  relaxation can potentially be used to fix variables at 0/1 values (see \cite{FLbook}, for example).  This is an important feature in the B{\&}B context. The methodology may be able to fix a number of variables when the relaxation generates a strong bound, and in doing so, it reduces the size of the  successive subproblems and improves the bounds computed for them.

In Table \ref{tab:fix}, we show the impact of using g-scaled linx, compared to o-scaled linx, on an iterative procedure where we solve \ref{linx}, \ref{DDFact}, and complementary \ref{DDFact}, fixing variables at 0/1 whenever possible. In both cases, we update the scaling parameter every time we solve \ref{linx}. For o-scaling, we optimize the scalar  $\gamma$ by applying Newton steps until  the absolute value of the derivative is less than $10^{-10}$. For g-scaling, we optimize the vector   $\Upsilon$ by applying up to 10 BFGS steps, taking $\gamma \mathbf{e}$ as a starting point. We limit the number of BFGS steps in this experiment to get closer to what might be practical within B\&B. We present in the columns of Table \ref{tab:fix}, the following information from left to right: The problem considered, $n$, the range of $s$ considered,
the scaling, the number of instances solved (one for each $s$ considered), the number of instances on which we could fix at least one variable (``inst fix''), the total number of variables fixed on all instances solved (``var fix''), the \%-improvement of g-scaling over o-scaling for the two last statistics. Additionally, to better understand how well our methods works for \hyperlink{MESP}{MESP} as $n$ grows, we also experimented with a covariance matrix of order $n=300$, which is a principal submatrix of the  covariance matrix of order $n=2000$ used as a benchmark in the literature (see \cite{li2020best,Fact}).   First, we see that, except for the number of instances of  \hyperlink{MESP}{MESP} with $n=124$ and $n=300$ on which we could fix variables,  there is always an improvement. The improvement becomes very significant when side constraints are considered.   We note that the number of variables fixed, reported on Table \ref{tab:fix}, refers only to the root nodes of the B{\&}B algorithm and indicates a promising approach to reduce  the B{\&}B enumeration.


\begin{table}[ht]
\setlength{\arraycolsep}{10pt}
\begin{center}
\begin{tabular}{l|ccc|ccc|rr}
&&&&\multicolumn{3}{c|}{Number of}&\multicolumn{2}{|c}{Improvement}\\
&$n$&s&scaling&$s$&inst fix&var fix&inst fix&var fix\\
\hline
\hyperlink{MESP}{MESP}
              & 63  & [2,62] & o & 61  & 41 & 1123   &          &          \\
              &     && g & 61  & 42 & 1140   & 2.44\%   & 1.51\%   \\
              & 90  & [2,89] & o & 88  & 41 & 1741   &          &          \\
              &     && g & 88  & 42 & 1790  & 2.44\%   & 2.81\%   \\
              & 124 & [2,123] & o & 122 & 35 & 3322 &          &          \\
              &     && g & 122 & 35 & 3353  & 0.00\%   & 0.93\%   \\
              & 300 & [80,120] & o & 41 & 41 &8382 &\\
              &     &&  g & 41 & 41 & 10753& 0.00\% &  28.3\%\\
\hline
\ref{CMESP}   & 63  & [3, 52] & o & 50  & 22 & 371    &          &          \\
              &     && g & 50  & 28 & 537    & 27.27\%  & 44.74\%  \\
              & 90  & [4, 87] & o & 84  & 26 & 606    &          &          \\
              &     && g & 84  & 37 & 1048   & 42.31\%  & 72.94\%  \\
              & 124 & [11, 110] & o & 100 & 9  & 197    &          &          \\
              &     && g & 100 & 33 & 1120  & 266.67\% & 468.53\%
\end{tabular}
\end{center}
\medskip
\caption{Impact of g-scaling on variable fixing}\label{tab:fix}
\end{table}

\vskip-10pt

\noindent The experiments with the fixing methodology
 show that  g-scaling  can effectively lead to a positive impact on the solution of \hyperlink{MESP}{MESP} and \ref{CMESP}, especially of the latter.

\section{Conclusion}\label{sec:conc}
We have seen that g-scaling can lead to improvements in upper bounds and variable fixing  for \hyperlink{MESP}{MESP} and very good improvements for \ref{CMESP}.
In future work, we will implement this in an efficient manner, within a B{\&}B algorithm.
In that context, it is important to efficiently use parent scaling vectors to warm-start
the optimization of scaling vectors for children (see \cite{Kurt_linx}). An open question
is whether g-scaling can help the \ref{DDFact} bound for \hyperlink{MESP}{MESP}.
Thm. \ref{thm:fact}.\emph{iv} is a partial result toward a negative answer.
Finally, there is another convex-optimization bound, the so-called ``NLP bound'' (see \cite{AFLW_Using}),
and it appears to be more difficult to get mathematical results on optimizing a g-scaling version of that
bound; but this is a good direction to explore.

\section{Appendix: Proof sketches}\label{sec:proofs}


\begin{proofsketch}[Thm. \ref{thm:bqp}]
\phantom{.}
\vspace{-7pt}

\begin{itemize}
    \item[\ref{thm:bqp}.i:] Suppose that the optimal solution to \ref{CMESP} is $x^*$, let $X^* := x^*\left(x^*\right)^\top$. Then $(x^*, X^*)\in P(n,S)$, and $f_{{\tiny\mbox{BQP}}}\left(x^*,X^*;\Upsilon\right) = \ldet C \left[S(x^*), S(x^*)\right]$. Thus $z_{{\tiny\mbox{BQP}}}(\Upsilon) \ge f_{{\tiny\mbox{BQP}}}\left(x^*,X^*;\Upsilon\right) = z$.
    \item[\ref{thm:bqp}.ii:] This is essentially a  result of \cite{Anstreicher_BQP_entropy}, with details filled in by \cite{FLbook}.
    \item[\ref{thm:bqp}.iii:] Let $F_{{\tiny\mbox{BQP}}}(x,X;\Upsilon):= \left(\Diag(\Upsilon)C\Diag(\Upsilon)\right)\circ X+ \Diag(\mathbf{e}-x)$ and $A_{{\tiny\mbox{BQP}}}(X;\Upsilon):=$\break  $ \left(\Diag(\Upsilon)C\Diag(\Upsilon)\right)\circ X$. Then given  $(x,X)$ in the domain of $f_{{\tiny\mbox{BQP}}}(x,X;\Upsilon)$ and feasible to \ref{BQP}, we have
     \begin{align*}
        &~\frac{\partial f^2_{{\tiny\mbox{BQP}}}\left(x,X;\Upsilon\right)}{\partial \left(\log \Upsilon\right)^2} \\
        =~&~ 4\Diag(x-\mathbf{e})\Diag\left(\diag\left(F_{{\tiny\mbox{BQP}}}(x,X;\Upsilon)^{-1}\right)\right)\\
        &~-4\Diag(x-\mathbf{e}) \left(F_{{\tiny\mbox{BQP}}}(x,X;\Upsilon)^{-1}\circ F_{{\tiny\mbox{BQP}}}(x,X;\Upsilon)^{-1}\right) \Diag(x-\mathbf{e}).
    \end{align*}
    When $x<\mathbf{e}$ and $X\succ 0$,  let $D_{{\tiny\mbox{BQP}}}(x):=\left(\Diag(\mathbf{e}-x)\right)^{1/2} \succ 0$ and further, $E_{{\tiny\mbox{BQP}}}(x,X;\Upsilon):=$\break
    $\left(D_{{\tiny\mbox{BQP}}}(x)\right)^{-1} A_{{\tiny\mbox{BQP}}}(X;\Upsilon) \left(D_{{\tiny\mbox{BQP}}}(x)\right)^{-1}\succ 0$. It can be shown that
    \begin{align*}
        &~\frac{\partial f^2_{{\tiny\mbox{BQP}}}\left(x,X;\Upsilon\right)}{\partial \left(\log \Upsilon\right)^2} = 4\left(E_{{\tiny\mbox{BQP}}}(x,X;\Upsilon)+I\right)^{-1}\circ \left(\left(E_{{\tiny\mbox{BQP}}}(x,X;\Upsilon)\right)^{-1}+I\right)^{-1}
        \succ  0.
    \end{align*}
    On the one hand, given $\Upsilon >0$, $\frac{\partial f^2_{{\tiny\mbox{BQP}}}\left(x,X;\Upsilon\right)}{\partial \left(\log \Upsilon\right)^2}$ is analytical on $(x,X)$ in the domain of $f_{{\tiny\mbox{BQP}}}\left(x,X;\Upsilon\right)$. On the other hand, the feasible set of \ref{BQP} is compact. Therefore, given $(x,X)$ in the domain of $f_{{\tiny\mbox{BQP}}}\left(x,X;\Upsilon\right)$ and feasible to \ref{BQP}, there exists $\epsilon>0$ such that $\mathcal{N}_{(x,X)}:= \left\{(x',X'): \|x-x'\|\le \epsilon\right\}$ $\cap \{ \text{domain of } f_{{\tiny\mbox{BQP}}}\left(x,X;\Upsilon\right)\}  \cap \{\text{feasible set to BQP}\} $ is compact. This implies that if $\frac{\partial f^2_{{\tiny\mbox{BQP}}}\left(x,X;\Upsilon\right)}{\partial \left(\log \Upsilon\right)^2}$ $\prec 0$, then $\exists\ (x',X')\in \mathcal{N}_{(x,X)}$ such that $x'<\mathbf{e}$ and $\frac{\partial f^2_{{\tiny\mbox{BQP}}}\left(x',X';\Upsilon\right)}{\partial \left(\log \Upsilon\right)^2}\prec 0$, a contradiction. So, for each fixed $(x,X)$ such above, $f_{{\tiny\mbox{BQP}}}\left(x,X;\Upsilon\right)$ is convex in $\log \Upsilon$. Because $z_{{\tiny\mbox{BQP}}}(\Upsilon)$ is the pointwise maximum over all $(x,X)\in P(n,x)$, it is convex in $\log \Upsilon$. \qed
\end{itemize}
\end{proofsketch}

\begin{proofsketch}[Thm. \ref{thm:linx}]
\phantom{.}
\vspace{-6pt}

\begin{itemize}
    \item[\ref{thm:linx}.i:] Suppose that the optimal solution to \ref{CMESP} is $x^*$; then we can show $f_{{\tiny\mbox{linx}}}(x^*;\Upsilon)  = $\break $ \ldet C \left[S(x^*), S(x^*)\right]$. Thus $z_{{\tiny\mbox{linx}}}(\Upsilon) \ge f_{{\tiny\mbox{linx}}}(x^*;\Upsilon)  = z$.
    \item[\ref{thm:linx}.ii:] This is essentially a  result of \cite{Kurt_linx}, with details filled in by \cite{FLbook}.
    \item[\ref{thm:linx}.iii:] Let $F_{{\tiny\mbox{linx}}}(x;\Upsilon):= \Diag(\Upsilon)C\Diag(x)C\Diag(\Upsilon)+ \Diag(\mathbf{e}-x)$ and $A_{{\tiny\mbox{linx}}}(x;\Upsilon):=$\break $ \Diag(\Upsilon)C\Diag(x)C\Diag(\Upsilon)$. Let $D_{{\tiny\mbox{linx}}}(x):=\left(\Diag(\mathbf{e}-x)\right)^{1/2}$ and $E_{{\tiny\mbox{linx}}}(x;\Upsilon)$\break  $:=$ $\left(D_{{\tiny\mbox{linx}}}(x)\right)^{-1} A_{{\tiny\mbox{linx}}}(x;\Upsilon) \left(D_{{\tiny\mbox{linx}}}(x)\right)^{-1}$  when $x< \mathbf{e}$. Then similar to \ref{thm:bqp}.iii. \qed

\end{itemize}
\end{proofsketch}

\begin{proofsketch}[Thm. \ref{thm:fact}]
\phantom{.}
\vspace{-5pt}

\begin{itemize}
    \item[\ref{thm:fact}.i:] This is essentially a  result of \cite{Fact}.
    \item[\ref{thm:fact}.ii:] This is essentially a  result of \cite{nikolov2015randomized}, with details filled in by \cite{FLbook}.
    \item[\ref{thm:fact}.iii:] Based on \cite[Proposition 2]{li2020best} and \cite[Theorem 2.4.18]{zalinescu2002convex}, we can show that for $x, \hat x$ in the domain of $f_{{\tiny\mbox{DDFact}}} (x;\Upsilon)$, the directional derivative of $f_{{\tiny\mbox{DDFact}}}(x;\Upsilon)$ at $x$ in direction $\frac{\hat x-x}{\|\hat x-x\|}$ is $T(x;\Upsilon)^\top \left( \frac{\hat x-x}{\|\hat x-x\|}\right)$ where
    \vspace{-2pt}
    \[
    T(x;\Upsilon):=\diag\left(F_{{\tiny\mbox{DDFact}}}(x;\Upsilon) Q \Diag\left( \beta ( \lambda )\right) Q ^\top F_{{\tiny\mbox{DDFact}}}(x;\Upsilon)^\top\right)-\log\Upsilon.
    \]

    We first show two preliminary results:
    \begin{itemize}
        \item[(a)] It can be shown that $f_{{\tiny\mbox{DDFact}}}(x;\Upsilon)$ is continuous on its domain. Then, because the feasible region of \ref{DDFact} is compact, given $x$, $\exists \tilde{r}>0$ such that $\forall r\le\tilde{r}$, $\mathcal{B}_r(x):=\{y: \|y-x\|\le r\}$ is included in the domain of $f_{{\tiny\mbox{DDFact}}}(x;\Upsilon)$. Furthermore, the intersection of the feasible region of \ref{DDFact} and $\mathcal{B}_r(x)$ is compact and included in the domain of $f_{{\tiny\mbox{DDFact}}}(x;\Upsilon)$, denoted as $\mathcal{N}_x^r$~, which implies  uniform continuity of $f_{{\tiny\mbox{DDFact}}}(x;\Upsilon)$ on $\mathcal{N}_x^r$~.
        \item[(b)] Let $\mathcal{C}(x):=\{y:\|y-x\|=1\}$. $\forall \epsilon>0$, by the Heine-Borel Theorem, $\exists$ a finite set $F\subset \mathcal{C}(x)$ such that $\forall y\in \mathcal{C}(x)$, $\exists
        u\in F$ such that $\|y-u\|<\epsilon$.
    \end{itemize}
    \smallskip
    Now we are ready to prove Thm. \ref{thm:fact}.iii. We will assume that $T(x;\Upsilon)\neq 0$ for simplicity. First, by the uniform continuity in (a), given $\epsilon>0$ and $r\le \tilde{r}$, $ \exists \delta\in (0, \epsilon)$ such that $\forall x_1, x_2 \in \mathcal{N}_x^r$ with $\|x_1-x_2\|\le \frac{\delta}{\|T(x;\Upsilon)\|}$, we have $|f_{{\tiny\mbox{DDFact}}}(x_1;\Upsilon)-f_{{\tiny\mbox{DDFact}}}(x_2;\Upsilon)|<\epsilon$.
 Second, by (b), $\exists F_{\epsilon}$ such that $\forall y\in \mathcal{C}(x)$, $\exists u\in \mathcal{C}(x)$ such that $\|y-u\|< \frac{\delta}{\|T(x;\Upsilon)\|\cdot \tilde{r}}$. Third, by the existence of directional derivatives of $f_{{\tiny\mbox{DDFact}}}(x;\Upsilon)$ at $x$, $\forall r\le\tilde{r}$ small enough, we have $\forall u\in F_{\epsilon}, t\le r$,
$
        \left|f_{{\tiny\mbox{DDFact}}}(x+tu;\Upsilon) - f_{{\tiny\mbox{DDFact}}}(x;\Upsilon)-t T(x;\Upsilon)^\top u\right|< \epsilon.
$
    Fourth, $\forall \hat x \in \mathcal{N}_x^r$, $\frac{\hat x}{\|\hat x\|}\in \mathcal{C}(x)$ and $\|\hat x\|\le r\le \tilde{r}$, and by the second argument, $\exists u\in F_{\epsilon}$ such that $\|\hat x- \|\hat x\|\cdot  u\|= \|\hat x\|\cdot \left\|\hat x/\|\hat x\|-u\right\|< \|\hat x\|\cdot \frac{\delta_1}{\|T(x;\Upsilon)\|\cdot \tilde{r}}\le \frac{\delta}{\|T(x;\Upsilon)\|}$.

    In all, given $\epsilon>0$, $\exists r\le \tilde{r}$ and $F_\epsilon$ such that $\forall \hat x\in \mathcal{N}_x^r$, $\exists u \in F_{\epsilon}$ such that
    \begin{align*}
         &\left|f_{{\tiny\mbox{DDFact}}}(\hat x;\Upsilon) - f_{{\tiny\mbox{DDFact}}}(x;\Upsilon)- T(x;\Upsilon)^\top (\hat x-x)\right|\\
         &\quad =~ \left|f_{{\tiny\mbox{DDFact}}}(\hat x;\Upsilon) - f_{{\tiny\mbox{DDFact}}}(\|\hat x\| \cdot u;\Upsilon)\right| \\
         &\qquad + \left|f_{{\tiny\mbox{DDFact}}}(\|\hat x\| \cdot u;\Upsilon) - f_{{\tiny\mbox{DDFact}}}(x;\Upsilon)- T(x;\Upsilon)^\top (\|\hat x\| \cdot u-x)\right|\\
         &\qquad +\left| T(x;\Upsilon)^\top(\|\hat x\| \cdot u-\hat x)\right|\\
         &\quad <~  \epsilon+\epsilon+\frac{\delta}{\|T(x;\Upsilon)\|}\cdot \|T(x;\Upsilon)\|
         ~<~   3\epsilon,
    \end{align*}
    which implies the result.
    \item[\ref{thm:fact}.iv:] By switching the role of $x$ and $\Upsilon$, we can show that for any $x$ in the domain of $f_{{\tiny\mbox{DDFact}}} (x;\Upsilon)$, there is a vector $\tilde{T}(x;\Upsilon)\in \mathbb{R}^n$ such that \begin{align*}
    \lim_{\|h\|\rightarrow 0 ~:~ \atop \Upsilon+h > 0} \frac{|f_{{\tiny\mbox{DDFact}}}(x;\Upsilon+h)-f_{{\tiny\mbox{DDFact}}}(x;\Upsilon)-\tilde{T}(x;\Upsilon)^\top h|}{\|h\|} = 0.
\end{align*}
When $\Upsilon>0$ falls into the interior of the positive cone, the above result is equivalent to  $f_{{\tiny\mbox{DDFact}}}(x;\mathbf{e})$ being differentiable in $\Upsilon$.

Letting $T(x^*;\Upsilon)$ be as defined in the proof of Thm. \ref{thm:fact}.iii, the remaining result is equivalent to $x^*\circ \left(T(x^*;\mathbf{e})-\mathbf{e}\right)=0$, which is further equivalent to
\[
\left(T(x^*;\mathbf{e})\right)_i=1, ~\forall x^*_i>0.
\]
Suppose that $\sigma$ is a permutation of $1,\cdots, n$ such that $\left(T(x^*;\mathbf{e})\right)_{\sigma(1)}\ge \cdots \ge \left(T(x^*;\mathbf{e})\right)_{\sigma(n)}$.
By \cite{li2020best} and KKT conditions for \ref{DDFact}, we have
\[
\sum_{i\in \{1,2,\ldots,n\}} x^*_{\sigma(i)} \left(T(x^*;\mathbf{e})\right)_{\sigma(i)}=\sum_{i\in  \{1,2,\ldots,s\}} \left(T(x^*;\mathbf{e})\right)_{\sigma(i)}=s.
\]
On the other hand, if $x^*_{\sigma(i)} = 1$, we have
\begin{align*}
    \left(T(x^*;\mathbf{e})\right)_{\sigma(i)} = ~&~F_{\sigma(i)\cdot} Q \Diag\left( \beta ( \lambda )\right) Q ^\top F_{\sigma(i)\cdot}^{\top} 
   ~\le~ 
    F_{\sigma(i)\cdot} \left(F_{{\tiny\mbox{DDFact}}}(x;\mathbf{e})\right)^{\dagger} F_{\sigma(i)\cdot}^{\top}\\
    = ~&~ \textstyle  F_{\sigma(i)\cdot} \left( F_{\sigma(i)\cdot}^\top F_{\sigma(i)\cdot}+ \sum_{j\neq \sigma(i)} x^*_j F_{j\cdot}^\top F_{j\cdot} \right)^{\dagger} F_{\sigma(i)\cdot}^{\top}\\
    \le ~&~ F_{\sigma(i)\cdot} \left( F_{\sigma(i)\cdot}^\top F_{\sigma(i)\cdot}\right)^{\dagger} F_{\sigma(i)\cdot}^{\top}
    ~=~ 1
\end{align*}
where the first inequality is due to $Q \Diag\left( \beta ( \lambda )\right) Q ^\top F_{{\tiny\mbox{DDFact}}}(x;\Upsilon)\preceq I$ and that the two matrices can be simultaneously  diagonalized by $Q$, and the second inequality is by the  Sherman–Morrison formula for the pseudo-inverse.

The above two formulae, together with the KKT conditions and $\sum_{i\in [n]} x^*_{\sigma(i)}=s$, imply that $\left(T(x^*;\mathbf{e})\right)_{\sigma(1)}= \cdots = \left(T(x^*;\mathbf{e})\right)_{\sigma(s)}=1$, and $\forall i >s$ such that $x^*_{\sigma(i)}>0$, $\left(T(x^*;\mathbf{e})\right)_{\sigma(i)}=\left(T(x^*;\mathbf{e})\right)_{\sigma(s)}=1$, which finishes the proof.
\qed

\end{itemize}
\end{proofsketch}






\subsubsection{Acknowledgements}
We are especially grateful to Kurt Anstreicher for suggesting the possibility of generalizing scaling for the linx bound.
\bibliographystyle{splncs04}
\bibliography{MO}
%




\end{document}